\documentstyle[amssymb,12pt]{amsart}%\input amssym.def

\newtheorem{prop}{Proposition}[section]

\newtheorem{lemma}{Lemma}[section]
\newtheorem{thm}{Theorem}[section]

\theoremstyle{remark}
\newtheorem{remark}{Remark}

\begin{document}
\newcommand{\nc}{\newcommand} \nc{\on}{\operatorname}
\nc{\pa}{\partial}
\nc{\cA}{{\cal A}}\nc{\cB}{{\cal B}}\nc{\cC}{{\cal C}}
\nc{\cE}{{\cal E}}\nc{\cG}{{\cal G}}\nc{\cH}{{\cal H}}
\nc{\cX}{{\cal X}}\nc{\cR}{{\cal R}}\nc{\cL}{{\cal L}}
\nc{\sh}{\on{sh}}\nc{\Id}{\on{Id}}
\nc{\ad}{\on{ad}}\nc{\Der}{\on{Der}}\nc{\End}{\on{End}}\nc{\res}{\on{res}}
\nc{\Imm}{\on{Im}}\nc{\limm}{\on{lim}}\nc{\Ad}{\on{Ad}}
\nc{\Hol}{\on{Hol}}\nc{\Det}{\on{Det}}
\nc{\de}{\delta}\nc{\si}{\sigma}\nc{\ve}{\varepsilon}
\nc{\al}{\alpha}
\nc{\CC}{{\Bbb C}}\nc{\ZZ}{{\Bbb Z}}\nc{\NN}{{\Bbb N}}
\nc{\AAA}{{\Bbb A}}\nc{\cO}{{\cal O}} \nc{\cF}{{\cal F}}
\nc{\la}{{\lambda}}\nc{\G}{{\frak g}}\nc{\A}{{\frak a}}
\nc{\HH}{{\frak h}}
\nc{\N}{{\frak n}}\nc{\B}{{\frak b}}
\nc{\La}{\Lambda}
\nc{\g}{\gamma}\nc{\eps}{\epsilon}\nc{\wt}{\widetilde}
\nc{\wh}{\widehat}
\nc{\bn}{\begin{equation}}\nc{\en}{\end{equation}}
\nc{\SL}{{\frak{sl}}}

% ****** GISPIC **********
%
%** by GISLI MASON *******
%
%**for commutative diagrams
%

\newcommand{\ldar}[1]{\begin{picture}(10,50)(-5,-25)
\put(0,25){\vector(0,-1){50}}
\put(5,0){\mbox{$#1$}} 
\end{picture}}

\newcommand{\lrar}[1]{\begin{picture}(50,10)(-25,-5)
\put(-25,0){\vector(1,0){50}}
\put(0,5){\makebox(0,0)[b]{\mbox{$#1$}}}
\end{picture}}

\newcommand{\luar}[1]{\begin{picture}(10,50)(-5,-25)
\put(0,-25){\vector(0,1){50}}
\put(5,0){\mbox{$#1$}}
\end{picture}}

\title[Hopf algebra cocycles for double Yangians]
{A construction of Hopf algebra
cocycles for the double Yangian $DY(\SL_{2})$}

\author{B. Enriquez}\address{Centre de Math\'ematiques, URA 169 du CNRS,
Ecole Polytechnique,
91128 Palai-seau, France}

\author{G. Felder}
\address{D-Math, ETH-Zentrum, HG G46, CH-8092 Z\"urich, Suisse}

\date{March 1997}
\maketitle

\begin{abstract} 
We construct a Hopf algebra cocycle in the Yangian double $DY(\SL_{2})$,
conjugating Drinfeld's coproduct to the usual one. To do that, we
factorize the twist between two ``opposite'' versions of Drinfeld's
coproduct, introduced in earlier work by V. Rubtsov and the first
author, using the decomposition of the algebra in its negative and
non-negative modes subalgebras. 
\end{abstract}

\subsection*{Introduction}

The purpose of this paper is to show that Drinfeld's coproduct of the
Yangian double $DY(\SL_{2})$ (\cite{D-new})
is conjugated to the usual one. 
For that, we construct a Hopf algebra cocycle in the Yangian double 
$DY(\SL_{2})$. 

Actually, we note that $DY(\SL_{2})$ is endowed with two variants of
Drinfeld's coproduct. These coproducts are associated with two
decompositions of the Lie algebra $\G = \SL_{2} \otimes \CC((z^{-1}))$, the 
first one being $\G = \G_{+} \oplus \G_{-}$, with $\G_{+} = (\HH \otimes
\CC[z]) \oplus \left( \N_{+}\otimes \CC((z^{-1})) \right)$, 
$\G_{-} = (\HH \otimes
z^{-1}\CC[[z^{-1}]]) \oplus \left( \N_{-}\otimes \CC((z^{-1})) \right)$,
and the
second one being its transform by the nontrivial Weyl group element of
$\SL_{2}$. 
Here $\HH$ and $\N_{\pm}$ are the standard Cartan and opposite
nilpotent subalgebras of $\SL_{2}$. 
In \cite{Enr-Rub}, we considered Hopf algebras $U_{\hbar}\G$
quantizing more general Lie
bialgebra structures associated with curves in higher genus, and showed
that they were conjugated by a twist $F$. 

The next step of \cite{Enr-Rub} was the construction of a deformation
$U_{\hbar} \G_{R}$ of
the enveloping algebra of an  algebra of regular functions with values
in $\SL_{2}$; in our ``rational'' situation, this Lie algebra
corresponds to $\SL_{2} \otimes \CC[z]$ and 
$U_{\hbar}\G_{R}$ to the  Yangian
$Y(\SL_{2})$. This subalgebra also had the property that 
$$
\Delta(U_{\hbar}\G_{R}) \subset U_{\hbar}\G \otimes U_{\hbar} \G_{R}, \quad
\bar\Delta(U_{\hbar}\G_{R}) \subset U_{\hbar}\G_{R} \otimes U_{\hbar} \G. 
$$
The last step of that paper was to decompose $F$ as a product
$$
F_{2}F_{1}, \quad \on{ with} \quad 
F_{1}\in U_{\hbar}\G \otimes U_{\hbar}\G_{R}, \quad 
F_{2}\in U_{\hbar}\G_{R} \otimes U_{\hbar}\G,
$$ 
and then to construct a
quasi-Hopf algebra structure on $U_{\hbar}\G_{R}$ by twisting the
coproduct $\Delta$ by $F_{1}$. $F_{1}$ and $F_{2}$ are constructed by
applying to one factor of $F$ a projection of $U_{\hbar}\G$ on
$U_{\hbar}\G_{R}$, which is a right $U_{\hbar}\G_{R}$-module map. 
In this construction, the choice of the projection is not
unique. Changing the projection has the effect of changing
$(F_{1},F_{2})$ into $(uF_{1},F_{2}u^{-1})$, for some $u\in
U_{\hbar}\G_{R}^{\otimes 2}$; this changes the coproduct $\Ad(F_{1})
\circ \Delta$ on $U_{\hbar}\G_{R}$ by some twist. 

The question naturally arises whether the same technique can be applied in
Hopf algebra situations. In this paper, we treat the case of the
rational Manin 
triple $\G = \G^{\ge 0}  \oplus \G^{<0}$, where $\G^{\ge 0} = \SL_{2}
\otimes \CC[z]$ and $\G^{< 0} = \SL_{2}
\otimes z^{-1}\CC[[z^{-1}]]$. 
In this situation,  both $\G^{\ge 0}$ and $\G^{<0}$ are Lie
subbialgebras of $\G$, and there are also deformations of their
enveloping algebras in $DY(\SL_{2})$, $A^{\ge 0} = Y(\SL_{2})$ and 
$A^{<0}$.  
Therefore we require that the projection $\Pi_{<0,r}$ 
be at the same time a right
$A^{<0}$-module map. Then it is uniquely determined. We show that
the the first part $F_{1}$ of the decomposition of $F$ 
contructed in this way satisfies the Hopf algebra cocycle
condition. This is the main result of this text. 

The proof of this fact
relies on the following results. We first prove that the second part
$F_{2}$ of the decomposition of $F$ is obtained by applying to $F$
a projection $\Pi_{\ge 0,l}$ 
similar to $\Pi_{<0,r}$ (eqs. (\ref{voron}), (\ref{nesher})). We
give two proofs (sections \ref{proof1},
\ref{proof2}) of
this result, both of them relying on some study of the
duality theory within $DY(\SL_{2})$ (section \ref{dualty}); the 
first proof directly applies results from \cite{Enr-Rub}. This enables
to show that the defect of the cocycle identity for $F_{1}$ belongs to
two spaces with intersection $1\otimes DY(\SL_{2}) \otimes 1$. The fact 
that the pentagon identity is automatically satisfied by such
defects (\cite{D-qH}) then shows that it is indeed equal to $1$. 

After we twist by $F_{1}$ the universal $R$-matrix
of $DY(\SL_{2})$ associated to Drinfeld's coproduct, we obtain 
a new solution of
the Yang-Baxter equation. Applying to it $2$-dimensional representations
of $DY(\SL_{2})$, we construct $L$-operators satisfying the Yangian
exchange (or $RLL$) relations of \cite{FRT,RS} (section \ref{sect:RLL}). 
This connection between Yangian $RLL$
relations and quantum current relations had earlier been obtained in
\cite{KhT} (see \cite{DF} in the trigonometric case). 
After this connection is clarified we are in position to 
show (section \ref{conjug}) that $F_{1}$ conjugates $\Delta$
to the Yangian coproduct on $DY(\SL_{2})$. 

We will consider an elliptic version of the construction of this paper
in a separate article (\cite{toap}). There we will construct `` twisted
cocycles'' providing solutions to the dynamical Yang-Baxter equation;
this will lead us to the construction of quantum currents of elliptic
quantum groups.

The first step towards generalizing our results to the case of a general
Lie algebra is to generalize the twist $F$. This has been done by
N. Reshetikhin (see Rem. \ref{atlas}). In the general case, $F$ is then
a product of factors coprersponding to each simple root; one might
expect that these factors satisfy braid relations. We would therefore
obtain a ``quantum currents'' version of braid group representations. 
The next step of that generalization would be the study of the duality
theory within general double Yangians.   

This work was done during our stay at the ``Semestre
syst\`emes int\'egrables'' organized at the Centre
Emile Borel, Paris, UMS 839, CNRS/UPMC. 
We would like to express our thanks to its organizers for their
invitation to this very stimulating meeting. We also would like to
acknowledge discussions with O. Babelon, D. Bernard and N. Reshetikhin.

\section{The double Yangian $DY(\frak{sl}_{2})$ and its coproducts}

The double Yangian is a Hopf algebra that was introduced in \cite{Kh}
(see also \cite{Sm} in the non centrally extended case), 
and is associated to any semisimple Lie algebra $\G$. In the case where 
$\G = \frak{sl}_{2}$, this algebra is denoted by 
$DY(\frak{sl}_{2})$. We will also denote it by $A$. 
It is an algebra over the ring of formal power series in the variable
$\hbar$, with generators $x_{n},n\in\ZZ$ ($x=e,f,h$), $D$ 
and $K$, and the following relations: 
\begin{equation} \label{k+:e}
k^{+}(z)e(w)k^{+}(z)^{-1} = {{z-w+\hbar}\over{z-w}}e(w),
\end{equation}
\begin{equation} \label{k+:f}
k^{+}(z)f(w)k^{+}(z)^{-1} = {{z-w}\over{z-w+\hbar}}f(w),
\end{equation}
\begin{equation} \label{k-:e}
k^{-}(z)e(w)k^{-}(z)^{-1} = {{w-z+\hbar K}\over{w-z}}e(w),
\end{equation}
\begin{equation} \label{k-:f}
k^{-}(z)f(w)k^{-}(z)^{-1} = {{z-w+\hbar}\over{z-w}}f(w),
\end{equation}
\begin{equation} \label{e:e}
(z-w - \hbar)e(z)e(w) = (z-w + \hbar) e(w) e(z), 
\end{equation}
\begin{equation} \label{f:f}
(z-w + \hbar)f(z)f(w) = (z-w - \hbar) f(w) f(z), 
\end{equation}
\begin{equation} \label{e:f}
[e(z),f(w)] = {1\over{\hbar}} \left(\delta(z,w)K^{+}(z) - 
\delta(z,w-\hbar K)
K^{-}(w)^{-1} \right) , 
\end{equation}
\begin{equation} \label{k:k}
[K,\on{anything}] = 0, \quad [k^{\pm}(z), k^{\pm}(w)] = 0, \quad
[D,x(z)] = {{dx}\over{dz}}(z), x=e,f,k^{\pm},
\end{equation}
\begin{align} \label{K:K}
(z-w-\hbar) & (z-w+\hbar-\hbar K) K^{+}(z) K^{-}(w) 
\\ & \nonumber = (z-w+\hbar)(z-w-\hbar+\hbar K) 
K^{-}(w) K^{+}(z), 
\end{align}
where for $x=e,f,h$, we set 
$$
x^{\ge 0}(z) = \sum_{n\ge 0} x_{n}z^{-n-1}, \quad 
x^{<0}(z) = \sum_{n < 0} x_{n}z^{-n-1}, \quad
x(z) = x^{\ge 0}(z) + x^{<0}(z);  
$$
we also set 
\begin{equation} \label{k+}
k^{+}(z) = \exp \left( h_{0}\ln({{z + \hbar}\over{z}})
+\sum_{n > 0} h_{n} {{z^{-n} - (z+\hbar)^{-n}}\over{n}}
\right) , 
\end{equation}
\begin{equation} \label{K-}
K^{-}(z) = \exp \left( \hbar \sum_{n < 0} h_{n} z^{-n-1} \right) , 
\end{equation}
and 
$K^{+}(z) = k^{+}(z)k^{+}(z-\hbar)$, $K^{-}(z) =
k^{-}(z)k^{-}(z-\hbar)$. In (\ref{k+}), (\ref{K-}), the arguments of the
exponentials are viewed as formal power series in $\hbar$, with
coefficients in 
$A\otimes z^{-1}\CC((z^{-1}))$ 
in the first case, and in 
$A\otimes \CC[z]$ in the second one. Finally, $\delta(z,w) =
\sum_{n\in\ZZ} z^{n}w^{-n-1}$. 

\begin{remark}
The $x_{n}$ correspond in the notation of \cite{Enr-Rub}, to $x[z^{n}]$, for
$x=e,f,h$ and $n\in\ZZ$. 
\hfill \qed \medskip\end{remark}

The double Yangian Hopf structure $\Delta_{Yg}$ is defined as follows 
(see \cite{Kh}). Set 
$$
L^{\ge 0}(z) = \pmatrix 1 & \hbar f^{\ge 0}(z) \\ 0 & 1 \endpmatrix 
\pmatrix
k^{+}(z-\hbar) & 0 \\ 0 & k^{+}(z)^{-1}\endpmatrix \pmatrix 1 & 0 \\
\hbar e^{\ge 0}(z) & 1 \endpmatrix
$$
and 
$$
L^{<0}(z) = \pmatrix 1 & 0 \\ \hbar e^{<0}(z-\hbar K) & 1 \endpmatrix
\pmatrix k^{-}(z-\hbar) & 0 \\ 0 & k^{-}(z)^{-1} \endpmatrix \pmatrix 1
& \hbar f^{<0}(z) \\ 0 & 1 \endpmatrix , 
$$
then $L^{\ge0, <0}(z)$ are formal series in $z$ with values in 
$A \otimes \End(\CC^{2})$, and we set 
$$
\Delta_{Yg}(K) = K \otimes 1 + 1 \otimes K, \quad
\Delta_{Yg}(D) = D \otimes 1 + 1 \otimes D, 
$$
\begin{equation} \label{Delta:Yg:1}
(\Delta_{Yg}\otimes 1)L^{\geq 0}(z) = L^{\geq 0}(z)^{(13)}
L^{\geq 0}(z)^{(23)}, 
\end{equation}
\begin{equation} \label{Delta:Yg:2}
(\Delta_{Yg}\otimes 1)L^{<0}(z) = L^{<0}(z-\hbar
K_{1})^{(23)}L^{<0}(z)^{(13)}, 
\end{equation}
with $K_{1} = K\otimes 1$. 

The algebra 
$A$ can also be endowed with Drinfeld's Hopf structures $(\Delta,
\varepsilon,S)$ and 
$(\bar\Delta,\varepsilon,\bar S)$. They are given, on the one hand,
by the coproduct $\Delta$ defined by 
\begin{equation} \label{Delta:k}
\Delta(k^{+}(z)) = k^{+}(z) \otimes k^{+}(z), \quad
\Delta(K^{-}(z)) = K^{-}(z) \otimes K^{-}(z +\hbar K_{1}), 
\end{equation}
\begin{equation} \label{Delta:e}
\Delta(e(z)) = e(z)\otimes K^{+}(z) + 1\otimes e(z),
\end{equation}
\begin{equation} \label{Delta:f}
\Delta(f(z)) = f(z)\otimes 1 + K^{-}(z)^{-1} \otimes f(z+\hbar K_{1}),
\end{equation}
\begin{equation} \label{Delta:D:K}
\Delta(D) = D \otimes 1 + 1\otimes D, \quad \Delta(K) = K \otimes 1 + 
1 \otimes K , 
\end{equation}
the counit $\varepsilon$, and the antipode $S$ defined by them; 
and on the other hand, by the coproduct $\bar\Delta$ defined by 
\begin{equation} \label{bar:Delta:k}
\bar\Delta(k^{+}(z)) = k^{+}(z) \otimes k^{+}(z), \quad
\bar\Delta(K^{-}(z)) = K^{-}(z) \otimes K^{-}(z + \hbar K_{1}), 
\end{equation}
\begin{equation} \label{bar:Delta:e}
\bar\Delta(e(z)) = e(z -\hbar K_{2})\otimes K^{-}(z -\hbar K_{2})^{-1} 
+ 1\otimes e(z),
\end{equation}
\begin{equation} \label{bar:Delta:f}
\bar\Delta(f(z)) = f(z)\otimes 1 + K^{+}(z) \otimes f(z),
\end{equation}
\begin{equation} \label{bar:Delta:D:K}
\bar\Delta(D) = D \otimes 1 + 1\otimes D, \quad \bar\Delta(K) = 
K \otimes 1 + 1 \otimes K , 
\end{equation}
the counit $\varepsilon$, and the antipode $\bar S$ defined by them. 

As we remarked in \cite{Enr-Rub}, $\Delta$ and $\bar\Delta$ are linked by a 
twist operation. Let us set 
$$
F = \exp\left( \hbar \sum_{n\in\ZZ} e_{n}\otimes f_{-n-1} \right) ,  
$$
then we have 
\begin{equation} \label{debdou}
\bar\Delta = \Ad(F) \circ \Delta. 
\end{equation}
(Here and later, we use the notation $\Ad(u)(x) = uxu^{-1}$, for
$x$ and $u$ elements of some algebra, with $u$ invertible.)

$F$ satisfies the cocycle condition
$$
(F\otimes 1)(\Delta\otimes 1)(F) = (1\otimes F)(1\otimes \Delta)(F). 
$$
(see \cite{Enr-Rub}). 

\begin{remark} \label{atlas}
N. Reshetikhin informed us that he obtained the
conjugation equation (\ref{debdou}) in the general case (that is, with 
$\SL_{2}$ replaced by a semisimple Lie algebra $\G$). Then $F$ is equal
to the product $\prod_{k=1}^{\nu}F_{i_{k}}$, where $w_{0}=
s_{i_{1}}\ldots s_{i_{\nu}}$ is a decomposition of the longest Weyl
group element as a product of simple reflections, and $F_{i} =
q^{\sum_{n\in\ZZ} e_{i;n} \otimes f_{i;-n-1}}$, with 
$(e_{i;n})_{n\in\ZZ},(f_{i;n})_{n\in\ZZ}$ the components of the fields
corresponding to the $i$th simple root. The
crossed vertex relations seem to imply that all elements $F_{i}$'s
commute together. However, this is not quite true: the relations 
$$
(z-w+\hbar a_{ij})(z-w-\hbar a_{ij}) \left[ 
e_{i}(z) \otimes f_{i}(z) , e_{j}(w) \otimes f_{j}(w) \right] 
= 0  
$$
(where $a_{ij}$ are the coefficients of the Cartan matrix)
do not imply that the fields $e_{i}(z)\otimes f_{i}(z)$ and
$e_{j}(w)\otimes f_{j}(w)$ commute, but rather the existence of fields 
$A_{ij}^{\pm}(z)$ such that 
$$
\left[ 
e_{i}(z) \otimes f_{i}(z) , e_{j}(w) \otimes f_{j}(w) \right] 
= \delta(z,w+\hbar a_{ij}) A_{ij}^{+}(z)+ 
\delta(z,w-\hbar a_{ij}) A_{ij}^{-}(z). 
$$
It would be interesting to check using these fields, whether the
$F_{i}$'s satisfy the braid relations. In the same spirit, one is led to
construct fields corresponding to non-simple roots using the relation
$(z-w-\hbar a_{ij})[e_{i}(z),(k_{j}e_{j})(w)]=0$. 
\end{remark}

\section{Decomposition of $F$}

\subsection{Subalgebras of $A$}

We will call $A^{\geq 0}$ and $A^{<0}$ the subalgebras of $A$ generated by 
$D$ and the $x_{n},n\ge 0$, resp. by $K$ and the $x_{n},n<0$ 
(with $x=e,f,h$.) The multiplication induces isomorphisms from
$A^{\geq 0}\otimes A^{<0}$ and $A^{<0} \otimes A^{\geq 0}$ 
to $A$; moreover, the intersection of $A^{\geq 0}$
with $A^{<0}$ is reduced to $\CC 1$. 

Let $U_{\hbar}\N_{+}$ and $U_{\hbar}\N_{-}$ be the subalgebras of $A$ 
generated by the $e_{n},n\in \ZZ$, resp. the $f_{n},n\in\ZZ$. 

Let $U_{\hbar}\N_{\eps}^{\ge 0}$ and $U_{\hbar}\N_{\eps}^{<0}$ be the
subalgebras generated by the $x_{n},n\ge 0$, resp. by the $x_{n},n<0$,
with $x = e$ for $\eps = +$ and $x = f$ for $\eps = -$. 

The linear maps $U_{\hbar}\N_{\eps}^{\ge 0} \otimes U_{\hbar}\N_{\eps}^{<0} 
\to U_{\hbar}\N_{\eps}$ and $U_{\hbar}\N_{\eps}^{<0} \otimes 
U_{\hbar}\N_{\eps}^{\ge 0} \to 
U_{\hbar}\N_{\eps}$, defined by the composition of the inclusion with
the multiplication, are linear isomorphisms; moreover, the inclusions of
algebras $U_{\hbar}\N_{\eps}^{\ge 0} \subset U_{\hbar}\N_{\eps}$ and
$U_{\hbar}\N_{\eps}^{< 0} \subset U_{\hbar}\N_{\eps}$  
are flat deformations of the inclusions of
commutative algebras $\CC[x_{n},n\ge 0] \subset \CC[x_{n},n\in\ZZ]$ and 
$\CC[x_{n},n <  0] \subset \CC[x_{n},n\in\ZZ]$ (see
e.g. \cite{Enr-Rub}). 

\begin{remark} Relations for generating currents $x^{\ge 0}(z) 
= \sum_{n\ge 0} x_{n}z^{-n-1}$ and $x^{<0}(z) 
= \sum_{n< 0} x_{n}z^{-n-1}$ of $U_{\hbar}\N_{\eps}^{\ge 0}$
and $U_{\hbar}\N_{\eps}^{< 0}$ are 
$$
(z-w - \hbar)e^{\eta}(z)e^{\eta}(w) 
- (z-w + \hbar)e^{\eta}(w)e^{\eta}(z) = - \hbar \left( e^{\eta}(z)^{2}
+e^{\eta}(w)^{2} \right) 
$$
and
$$
(z-w + \hbar)f^{\eta}(z)f^{\eta}(w) 
- (z-w - \hbar)f^{\eta}(w)f^{\eta}(z) =  \hbar \left( f^{\eta}(z)^{2}
+f^{\eta}(w)^{2} \right) , 
$$
$\eta\in\{\ge 0, <0\}$. 

On the other hand, the relations between these currents are 
$$ 
(z-w-\hbar) e^{\eta}(z) e^{\eta'}(w)
-
(z-w+\hbar) e^{\eta'}(w) e^{\eta}(z)
-
\hbar[e^{\eta}(z)^{2} + e^{\eta'}(w)^{2}] = 0 ,
$$ 
$$ 
(z-w+\hbar) f^{\eta}(z) f^{\eta'}(w)
-
(z-w-\hbar) f^{\eta'}(w) f^{\eta}(z)
+
\hbar[f^{\eta}(z)^{2} + f^{\eta'}(w)^{2}] = 0 . 
$$ 
if $\{\eta,\eta'\} = \{\ge 0, < 0\}$. 
\end{remark}

\subsection{Hopf algebra pairings} \label{dualty}

Let $U_{\hbar}\HH_{+}$ be the subalgebra of $A$ generated by $D$ and the
$h_{n},n\ge 0$, and $U_{\hbar}\HH_{-}$ be the subalgebra of $A$
generated by $K$ and the $h_{n},n< 0$,

Let $U_{\hbar}\G_{\pm}$ be the subalgebras of $A$ generated by
$U_{\hbar}\HH_{\pm}$ and $U_{\hbar}\N_{\pm}$, and $U_{\hbar}\bar\G_{\pm}$
the subalgebras of $A$ generated by $U_{\hbar}\HH_{\mp}$ and 
$U_{\hbar}\N_{\pm}$. 

$(U_{\hbar}\G_{\pm},\Delta)$ are Hopf subalgebras of $(A,\Delta)$;
$(U_{\hbar}\G_{+},\Delta)$ and $(U_{\hbar}\G_{-},\Delta')$
are dual to each other, and the duality $\langle , \rangle$ is 
expressed by the rules 
$$
\langle e_{n} , f_{m} \rangle = {1\over \hbar}\delta_{n+m+1,0}, \quad
\langle h_{a} , h_{b} \rangle = {2\over \hbar}\delta_{a+b+1,0}, \quad
\langle D,K\rangle = {1\over\hbar}, 
$$
$n,m\in\ZZ,a\ge 0,b < 0$, 
the other pairings between generators being trivial. 

In a similar way, $(U_{\hbar}\bar\G_{\pm},\bar\Delta)$ are Hopf 
subalgebras of $(A,\Delta)$;
$(U_{\hbar}\bar\G_{+},\bar\Delta')$ and $(U_{\hbar}\bar\G_{-},$
$\bar\Delta)$
%%%%%%%%%%%% LIAISON
are dual to each other, and the duality $\langle , \rangle'$ is 
expressed by the rules 
$$
\langle e_{n} , f_{m} \rangle' = {1\over \hbar}\delta_{n+m+1,0}, \quad
\langle h_{a} , h_{b} \rangle' = {2\over \hbar}\delta_{a+b+1,0}, \quad
\langle K,D\rangle' = {1\over\hbar}, 
$$
$n,m\in\ZZ,a\ge 0,b < 0$, 
the other pairings between generators being trivial. 

The restrictions of $\langle , \rangle$ and $\langle , \rangle'$ to 
$U_{\hbar}\N_{+}\times U_{\hbar}\N_{-}$  coincide and are denoted by 
$\langle , \rangle_{U_{\hbar}\N_{\pm}}$.

Moreover, we have

\begin{lemma} \label{orthog} (see \cite{Enr-Rub})
1) The annihilator of $U_{\hbar}\N_{-}^{\ge 0}$ for $\langle ,
\rangle_{U_{\hbar}\N_{\pm}}$ is $\sum_{n\ge 0}e_{n}
 \cdot U_{\hbar}\N_{+}$.

2) The annihilator of $U_{\hbar}\N_{+}^{<0}$ for $\langle ,
\rangle_{U_{\hbar}\N_{\pm}}$ is $\sum_{n< 0}f_{n}\cdot U_{\hbar}\N_{-}$.

3) The annihilator of $U_{\hbar}\N_{+}^{\ge 0}$ is $\sum_{n\ge 0}
U_{\hbar}\N_{-} \cdot f_{n}$.

4) The annihilator of $U_{\hbar}\N_{-}^{<0}$ for $\langle ,
\rangle_{U_{\hbar}\N_{\pm}}$ is $\sum_{n<0}U_{\hbar}\N_{+}\cdot e_{n}$.
\end{lemma}

{\em Proof.\/} 1) and 3) are consequences of \cite{Enr-Rub}, Prop. 6.2, 
and 2) and 4) are shown in a similar way.
\hfill \qed \medskip 

Finally, the link between $F$ the pairing $\langle ,
\rangle_{U_{\hbar}\N_{\pm}}$ can be described as follows. Let us first
introduce the notation
$$
\langle a, id \otimes b \rangle_{V,W} = 
\sum_{i}a_{i}\langle a'_{i},b\rangle_{V,W},  \quad
\langle a, b \otimes id \rangle_{V,W} = 
\sum_{i}\langle a_{i},b\rangle_{V,W} a'_{i},  
$$
for $a \in V^{\otimes 2}$ and $b\in W$, for $V,W$ some vector spaces and
$\langle , \rangle_{V,W}$ some pairing between them, $a$ being
decomposed as $\sum_{i}a_{i}\otimes a'_{i}$.

\begin{lemma} \label{F-pair}
(see \cite{Enr-Rub}, (66) and (68))
1) For any $x\in U_{\hbar}\N_{+}$, we have 
$$
\langle F , id \otimes x \rangle_{U_{\hbar}\N_{\pm}} = x. 
$$
2) For any $y\in U_{\hbar}\N_{-}$, we have 
$$
\langle F , y \otimes id \rangle_{U_{\hbar}\N_{\pm}} = y. 
$$ 
\end{lemma}

\subsection{Decomposition of $F$} \label{proof1}

\begin{prop} \label{aryeh}
There exists a decomposition $F = F_{2} F_{1}$, with 
$F_{1}\in U_{\hbar}\N_{+}^{<0}
\otimes U_{\hbar}\N_{-}^{\ge 0}
$ and 
$F_{2} \in U_{\hbar}\N_{+}^{\ge 0}
\otimes U_{\hbar}\N_{-}^{< 0}$. 
It is unique up to changes of
$(F_{1},F_{2})$ into $(\la F_{1},\la^{-1}F_{2})$, with
$\la\in\CC^{\times}$.  
\end{prop}

{\em Proof. \/}
Let us denote by 
$\Pi_{\ge 0,l}$,  $\Pi_{\ge 0,r}$; and by   
$\Pi_{<0,l}$,  $\Pi_{<0,r}$
the linear maps from $U_{\hbar}\N_{\eps}$ to
$U_{\hbar}\N_{\eps}^{\ge 0}$ and $U_{\hbar}\N_{\eps}^{<0}$ 
defined by
$$
\Pi_{\eta,l}(a_{\eta}a_{\eta'}) = a_{\eta}\varepsilon(a_{\eta'}), 
\quad
\Pi_{\eta,r}(a_{\eta'}a_{\eta}) = \varepsilon(a_{\eta'})a_{\eta},  
$$
for $\{\eta,\eta'\} = \{\ge 0, <0 \}$ and 
$a_{\eta}\in U_{\hbar}\N_{\eps}^{\eta}$.

\begin{lemma}
1)  $(\Pi_{<0,r}\otimes 1)(F)$ belongs to $U_{\hbar}\N_{+}^{<0}
\otimes U_{\hbar}\N_{-}^{\ge 0}$. 

2) $(1\otimes \Pi_{\ge 0,r})(F)$ belongs to 
$U_{\hbar}\N_{+}^{<0}\otimes U_{\hbar}\N_{-}^{\ge 0}$. 
\end{lemma}

{\em Proof. \/}
1) $(\Pi_{<0,r}\otimes 1)(F)$ clearly  
belongs to $U_{\hbar}\N_{+}^{<0} \otimes U_{\hbar}\N_{-}$. 
On the other hand, 
we have for any $a\in U\N_{+}$ and $n\ge 0$, 
$$
\langle (\Pi_{<0,r}\otimes 1)(F), id\otimes e_{n} a\rangle = 
\Pi_{<0,r}(e_{n} a) = 0 ; 
$$
the first equality follows from Lemma \ref{F-pair}, 1, and the second
from the fact that $\Pi_{<0,r}$ is a left $U_{\hbar}\N_{+}^{\ge
0}$-module map. 
From Lemma \ref{orthog}, 1,  now follows that 
 $(\Pi_{<0,r}\otimes 1)(F)$ also belongs to 
$U_{\hbar}\N_{+}\otimes U_{\hbar}\N_{-}^{\ge 0}$.

2) is proved in the same way, using Lemma \ref{F-pair}, 2, and Lemma
\ref{orthog}, 2. 
\hfill \qed \medskip

\begin{lemma}
$(\Pi_{<0,r}\otimes 1)(F)$ is equal to 
$(1\otimes \Pi_{\ge 0,r})(F)$. 
\end{lemma}

{\em Proof. \/}
Let $a_{+}$ belong to $U_{\hbar}\N_{-}^{\ge 0}$ 
and let $a_{-}$ belong to $U_{\hbar}\N_{+}^{< 0}$. Let us compute 
\begin{equation} \label{pair***}
\langle (\Pi_{<0,r}\otimes 1)(F) - (1\otimes \Pi_{\ge 0,r})(F) , 
a_{+}\otimes a_{-}\rangle_{U_{\hbar}\N_{\pm}^{\otimes 2}}.
\end{equation} 

Due to Lemma \ref{F-pair}, this is equal to $\langle \Pi_{<0,r}(a_{-}), 
a_{+}\rangle_{U_{\hbar}\N_{\pm}} - \langle a_{-}, \Pi_{\ge 0,r}(a_{+})
\rangle_{U_{\hbar}\N_{\pm}}$. 

Since $\Pi_{<0,r}(a_{-}) = a_{-}$, $\Pi_{\ge 0,r}(a_{+}) =
a_{+}$, (\ref{pair***}) is equal to zero. 

The pairing $\langle , \rangle_{U_{\hbar}\N_{\pm}}$ is a flat 
deformation of the symmetric power of the pairing between $\CC[z]$ and
$z^{-1}\CC[[z^{-1}]]$, defined by $\langle f,g \rangle =
\res_{\infty}(fg dz)$. Therefore, it defines an
injection of $U_{\hbar}\N_{-}^{\ge 0}$ in the dual of
$U_{\hbar}\N_{+}^{<0}$ and of $U_{\hbar}\N_{+}^{<0}$ in the dual of 
$U_{\hbar}\N_{-}^{\ge 0}$. That (\ref{pair***}) is equal to zero then
implies that 
$(\Pi_{<0,r}\otimes 1)(F) = (1\otimes \Pi_{\ge 0,r})(F)$. 
\hfill \qed \medskip

As $\Pi_{\ge 0,r}$ is a right $U_{\hbar}\N_{+}^{\ge 0}$-module map, we can 
apply \cite{Enr-Rub}, (74), second statement, with $U_{\hbar}\G_{R} =
A^{\geq 0}$, 
and obtain  % jeru 
\begin{equation} \label{propr1}
(1\otimes \Pi_{\ge 0 ,r})(F) F^{-1} \in U_{\hbar}\N_{+}^{\ge 0} 
\otimes U_{\hbar}\N_{-}.
\end{equation}
We can now apply the arguments of \cite{Enr-Rub}, Prop. 7.2, to the Hopf 
algebra $(A,\Delta')$. The role of $U_{\hbar}\G_{R}$ is now played 
by $A^{<0}$; $F$ is replaced by $F^{(21)}$. The analogue of the second 
statement of \cite{Enr-Rub}, (74) is then 
\begin{equation} \label{propr2}
( (1\otimes \Pi_{<0,r})(F^{(21)})) (F^{(21)})^{-1} \in
U_{\hbar}\N_{-}^{<0} \otimes U_{\hbar}\N_{+}. 
\end{equation}
We can show in a similar way that 
$$
(\Pi_{\ge 0,l}\otimes 1)(F) = (1\otimes \Pi_{<0,l})(F),
$$
so that this quantity belongs to $U_{\hbar}\N_{+}^{\ge 0}\otimes 
U_{\hbar}^{<0}\N_{-}$, and that
\begin{equation} \label{propl}
F^{-1} (\Pi_{\ge 0,l}\otimes 1)(F) \in U_{\hbar}\N_{+}^{< 0} \otimes 
U_{\hbar}\N_{-}^{\ge 0}. 
\end{equation}

Consider now the product 
\begin{equation} \label{prod}
((\Pi_{\ge 0,l}\otimes 1)(F))^{-1} F 
((\Pi_{<0,r}\otimes 1)(F)
)^{-1}.
\end{equation}
Since  
$(\Pi_{<0,r}\otimes 1)(F) \in U_{\hbar}\N_{+}^{<0} \otimes 
U_{\hbar}\N_{-}^{\ge 0}$, and by  (\ref{propl}), this product belongs to 
$U_{\hbar}\N_{+}^{<0} \otimes U_{\hbar}\N_{-}^{\ge 0}
$. 
On the other hand, since 
$(\Pi_{\ge 0,l}\otimes 1)(F)\in U_{\hbar}\N_{+}^{\ge 0} 
\otimes U_{\hbar}\N_{-}^{<0}$, and by (\ref{propr1}), 
it belongs to 
$U_{\hbar}\N_{+}^{<0} \otimes U_{\hbar}\N_{-}^{\ge 0}$. It follows 
that this product is scalar. Since the constant term in its expansion is
equal to one, (\ref{prod}) is equal to one. 

Therefore we can set 
\begin{equation} \label{voron}
F_{1} = (\Pi_{<0,r}\otimes 1)(F) = (1\otimes \Pi_{\ge 0,r})(F) 
\end{equation}
and 
\begin{equation} \label{nesher}
F_{2} = (\Pi_{\ge 0,l}\otimes 1)(F) = (1\otimes \Pi_{<0,l})(F). 
\end{equation}
\hfill \qed\medskip 

\subsection{Another proof of Prop. \ref{aryeh}.} \label{proof2}
Let us define $F_{1}$ and $F_{2}$ by 
\begin{equation} \label{F1:2}
F_{1} = (\Pi_{<0,r}\otimes 1)(F), \quad  F_{2} = 
(\Pi_{\ge 0,l}\otimes 1)(F), 
\end{equation}
and show directly that $F = F_{2} F_{1}$. For this, we will consider
the linear endomorphism $\ell$ of $U_{\hbar}\N_{+}$ defined by 
\begin{equation} \label{ell}
\ell(x) = \langle F_{2}F_{1}, id \otimes x\rangle_{U_{\hbar}\N_{\pm}}. 
\end{equation}

Let us denote by $\pi$ the linear map from $U_{\hbar}\G_{+}$ to 
$U_{\hbar}\N_{+}$, defined by $\pi(tx) = \varepsilon(t) x$, for 
$x\in U_{\hbar}\N_{+}$, $t\in U_{\hbar}\HH_{+}$. Let us also denote
by $\pi'$ the linear map from $U_{\hbar}\bar\G_{+}$ to 
$U_{\hbar}\N_{+}$, defined by $\pi'(x't') = x'\varepsilon(t')$, for 
$x'\in U_{\hbar}\N_{+}$, $t'\in U_{\hbar}\HH_{-}$. 

\begin{lemma} \label{gogol}
1) For $y\in U_{\hbar}\G_{+}$, we have
\begin{equation} \label{sfax}
\langle F, id \otimes y \rangle = \pi(y).  
\end{equation}

2) For $z\in U_{\hbar}\bar\G_{+}$, we have
\begin{equation} \label{cassini}
\langle F, id \otimes z \rangle' = \pi'(z).  
\end{equation}
\end{lemma}

{\em Proof. \/}
Let us prove 1). 
Let us first show that 
for any $y'\in U_{\hbar}\N_{-}$, we have 
\begin{equation} \label{piombi}
\langle y', y\rangle  = \langle y' , \pi(y) \rangle.
\end{equation} 
To prove this, consider the case where $y = t_{0} y_{0}$, 
$y_{0}\in U_{\hbar}\N_{+}$, $t_{0}\in U_{\hbar}\HH_{+}$. Then 
$\langle y', y\rangle  = \langle \Delta'(y') , 
t_{0} \otimes y_{0} \rangle_{(2)}$; but $\Delta'(y')$ belongs to 
$U_{\hbar}\N_{-} \otimes U_{\hbar}\G_{-}$, and for 
$a\in U_{\hbar}\HH_{+}$, $b\in U_{\hbar} \N_{-}$, $\langle a,b \rangle 
= \varepsilon(a) \varepsilon(b)$. 
It follows that  
$\langle y', y\rangle  = \langle (\varepsilon\otimes 1) \circ \Delta'(y')
\varepsilon(t_{0}), y_{0} \rangle = \langle y' , \pi(y) \rangle$, so that 
(\ref{piombi}) holds. (\ref{sfax}) then follows from (\ref{piombi}) and 
Lemma \ref{F-pair}. 

2) is proved in a similar way. 
\hfill \qed \medskip

We then compute $\ell(x)$ as follows, 
for $x\in U_{\hbar}\N_{+}$. Set $\Delta(x) = \sum_{i}x'_{i} \otimes
x''_{i}$, with $x'_{i}\in U_{\hbar}\N_{+}, x''_{i}\in U_{\hbar}\G_{+}$. 
Then
\begin{align} \label{expr:ell} \nonumber
& \ell(x) = \sum_{i} 
\langle F_{2}, id \otimes x'_{i} \rangle 
\langle F_{1}, id \otimes x''_{i} \rangle
= \sum_{i} \Pi_{\ge 0, l}(x'_{i}) \Pi_{<0,r}(\langle F, id \otimes x''_{i}
\rangle)
\\ & = \sum_{i} \Pi_{\ge 0, l}(x'_{i}) (\Pi_{<0,r} \circ \pi)(x''_{i}) . 
\end{align}
We deduce from this expression the following property of $\ell$. 

\begin{lemma} \label{tinto}
$\ell$ is a left $U_{\hbar}\N_{+}^{\ge 0}$-module map. 
\end{lemma}

{\em Proof. \/}
$\Pi_{<0,r}\circ \pi$ is defined as follows. 
Recall that the product operation defines a linear isomorphism from 
the tensor
product $U_{\hbar}\HH_{+} \otimes U_{\hbar}\N_{+}^{\ge 0} \otimes 
U_{\hbar}\N_{+}^{<0}$ onto $U_{\hbar}\G_{+}$. 
 $\Pi_{<0,r}\circ \pi$ is then defined by 
$(\Pi_{<0,r}\circ \pi)(x)
= \varepsilon(tx_{\ge 0})x_{<0}$, for $x$ decomposed as 
 $tx_{\ge 0}x_{<0}$, $t\in U_{\hbar}\HH_{+}$, 
$x_{\ge 0}\in U_{\hbar}\N_{+}^{\ge 0}, x_{<0}\in U_{\hbar}\N_{+}^{<0}$.
On the other hand, denote by $U_{\hbar}\B_{+}$ the subspace of 
$U_{\hbar}\G_{+}$ corresponding to 
$U_{\hbar}\HH_{+} \otimes U_{\hbar}\N_{+}^{\ge 0} \otimes 1$. 
We can check that this is a subalgebra of $U_{\hbar}\G_{+}$. 
It follows that $\Pi_{<0,r}\circ \pi$ satisfies
\begin{equation} \label{module}
(\Pi_{<0,r}\circ\pi) (bx) = \varepsilon(b) (\Pi_{<0,r}\circ\pi)(x), 
\end{equation}
for $b\in U_{\hbar}\B_{+}$, $x\in U_{\hbar}\G_{+}$. 

Finally, (\ref{Delta:e}) implies that for any $n\in\ZZ$, 
$\Delta(e_{n}) =  1\otimes e_{n} + \sum_{p\ge 0} e_{n-p} \otimes K^{+}_{p}$
(we set $K^{+}(z) = \sum_{p\ge 0} K^{+}_{p}z^{-p}$),
so that for $n\ge 0$ this belongs to $U_{\hbar}\N_{+} \otimes 
U_{\hbar}\B_{+}$; it follows that 
$\Delta(U_{\hbar}\N_{+}^{\ge 0})
\subset U_{\hbar}\N_{+}\otimes U_{\hbar}\B_{+}$.

Let us fix then $b$ in $U_{\hbar}\N_{+}^{\ge 0}$ and $x$ in 
$U_{\hbar}\N_{+}$. 
Set $\Delta(x) = \sum_{i} x'_{i} \otimes x''_{i}$, $x'_{i} \in U_{\hbar}
\N_{+}, x''_{i}\in U_{\hbar}\G_{+}$, and 
$\Delta(b) = \sum_{j} b'_{j} \otimes b''_{j}$, $b'_{j} \in U_{\hbar}
\N_{+}, b''_{j}\in U_{\hbar}\B_{+}$. 
Then
\begin{align*}
\ell(bx) & = 
\sum_{i,j} \Pi_{\ge 0, l}(b'_{j}x'_{i}) 
(\Pi_{<0,r} \circ \pi)(b''_{j}x''_{i}) = 
\sum_{i,j} \Pi_{\ge 0, l}(b'_{j}x'_{i}) 
\varepsilon(b''_{j}) (\Pi_{<0,r} \circ \pi) (x''_{i})
\\ & 
= \sum_{i} \Pi_{\ge 0, l}( b x'_{i}) 
(\Pi_{<0,r} \circ \pi)(x''_{i}) = 
b \sum_{i} \Pi_{\ge 0, l}( x'_{i}) 
(\Pi_{<0,r} \circ \pi)(x''_{i}) \\ & 
= 
b\ell(x) ;  
\end{align*}
the second equality follows from (\ref{module}), the third from 
the properties of the counit, and the fourth from the  fact that 
$\Pi_{\ge 0,r}$ is a left 
$U_{\hbar}\N_{+}^{\ge 0}$-module map. 
\hfill \qed \medskip

Set now $\bar\Delta(x) = \sum_{i}\bar x'_{i} \otimes
\bar x''_{i}$, with $\bar x'_{i}\in U_{\hbar}\N_{+}, 
\bar x''_{i}\in U_{\hbar}\bar \G_{+}$. 
Then
\begin{align} \label{expr:ell:2} \nonumber
\ell(x) & = \sum_{i} 
\langle F_{2}, id\otimes \bar x''_{i} \rangle' 
\langle F_{1}, id\otimes \bar x'_{i} \rangle'
= \sum_{i} \Pi_{\ge 0, l}( \langle F, id \otimes \bar x''_{i} \rangle') 
\Pi_{<0,r} (\bar x'_{i})
\\ & = \sum_{i} (\Pi_{\ge 0, l} \circ\pi')(\bar x''_{i}) 
\Pi_{<0,r}(\bar x'_{i}). 
\end{align}
We now deduce from this expression: 

\begin{lemma} \label{retto}
$\ell$ is a right $U_{\hbar}\N_{+}^{< 0}$-module map. 
\end{lemma}

{\em Proof. \/}
As above, the product operation defines an isomorphism of vector spaces 
from 
$U_{\hbar}\N_{+}^{\ge 0} \otimes U_{\hbar}\N_{+}^{<0} \otimes 
U_{\hbar}\HH_{-}$ 
to $U_{\hbar}\bar\G_{+}$. 
The image by this map of 
$1 \otimes U_{\hbar}\N_{+}^{<0} \otimes U_{\hbar}\HH_{-}$ is a subalgebra
of $U_{\hbar}\bar\G_{+}$, that we denote by $U_{\hbar}\bar\B_{+}$. 
$\Pi_{\ge 0, l} \circ\pi'$ is then defined by 
$(\Pi_{\ge 0, l} \circ\pi')(x) = \sum_{\al}x_{>0 ; \al} 
\varepsilon(b_{\al})$, 
if $x$ is decomposed as $\sum_{\al}   x_{>0 ; \al}b_{\al}$, 
$x_{>0;\al} \in U_{\hbar}\N_{+}^{\ge 0}$, $b_{\al}\in U_{\hbar}
\bar \B_{+}$. 
Therefore, we have
\begin{equation} \label{module2}
(\Pi_{\ge 0, l} \circ\pi')(xb) = (\Pi_{\ge 0, l} \circ\pi')(x)
\varepsilon(b), 
\end{equation}
for $x\in U_{\hbar}\bar\G_{+}, b\in U_{\hbar}\bar\B_{+}$. 

Finally, (\ref{bar:Delta:e}) implies that for $n>0$, 
$\bar\Delta(e_{-n}) = \sum_{p\ge 0} e_{-n+p} \otimes ((K^{-})^{-1})_{-p}
+ 1\otimes e_{-n}$ (we set $(K^{-})^{-1}(z) = \sum_{p\leq 0}
((K^{-})^{-1})_{p}z^{-p}$), and so belongs to $U_{\hbar}\N_{+} \otimes 
U_{\hbar}\bar \B_{+}$; it follows that 
$\bar\Delta(U_{\hbar}\N_{+}) \subset U_{\hbar}\N_{+}
\otimes U_{\hbar}\bar\B_{+}$. 

Fix then $x$ in $U_{\hbar}\N_{+}$, $b$ in $U_{\hbar}\N_{+}^{<0}$, 
with $\bar\Delta(x) = \sum_{i}\bar x'_{i}\otimes x''_{i}$, 
$x'_{i} \in U_{\hbar} \N_{+}$, $x''_{i} \in U_{\hbar}\bar\G_{+}$,
$b'_{i} \in U_{\hbar} \N_{+}^{< 0}$, $b''_{i} \in U_{\hbar}\bar\B_{+}$,
$\bar\Delta(b) = \sum_{j} \bar b'_{j} \otimes \bar b''_{j}$.
Then 
\begin{align*}
\ell(xb) & = \sum_{i,j} (\Pi_{\ge 0, l} \circ\pi')(\bar x''_{i}\bar 
b''_{j}) 
\Pi_{<0,r}(\bar x'_{i}\bar b'_{j})
= 
\sum_{i,j} (\Pi_{\ge 0, l} \circ\pi')(\bar x''_{i})\varepsilon
(\bar b''_{j}) 
\Pi_{<0,r}(\bar x'_{i}\bar b'_{j})
\\ & =
\sum_{i,j} (\Pi_{\ge 0, l} \circ\pi')(\bar x''_{i}) 
\Pi_{<0,r}(\bar x'_{i}b)
=
\sum_{i} (\Pi_{\ge 0, l} \circ\pi')(\bar x''_{i})
\Pi_{<0,r}(\bar x'_{i})b
\\ & 
= \ell(x) b; 
\end{align*}
the second equality follows from (\ref{module2}), the third one from 
the properties of $\varepsilon$, and the fourth one from 
the fact that $\Pi_{<0,r}$ is a right $U_{\hbar}\N_{+}^{<0}$-module map. 
\hfill \qed \medskip

Let us now prove Prop. \ref{aryeh}. 
We have $\ell(1) = 1$. 
Since any element of $U_{\hbar}\N_{+}$ can be 
expressed as a sum of products $\sum_{i}x_{i}^{\ge 0}x_{i}^{<0}$, 
with $x_{i}^{\ge 0} \in U_{\hbar}\N_{+}^{\ge 0}$, 
$x_{i}^{<0} \in U_{\hbar}\N_{+}^{<0}$, 
and by Lemmas \ref{tinto} and
\ref{retto}, $\ell$ coincides with the identity. \hfill \qed

\section{Cocycle properties}

\begin{thm} \label{cocycle}  
$F_{1}$ satisfies the cocycle equation
$$
(F_{1}\otimes 1)(\Delta\otimes 1)(F_{1}) = (1\otimes F_{1})(1\otimes
\Delta)(F_{1}).  
$$
\end{thm}

{\em Proof. \/}
First note that 
$$
\Delta(A^{\geq 0}) \subset A \otimes A^{\geq 0}, \quad
\Delta(A^{<0}) \subset A^{<0} \otimes A, 
$$
$$
\bar \Delta(A^{\geq 0}) \subset A^{\geq 0} \otimes A, \quad
\bar \Delta(A^{<0}) \subset A \otimes A^{<0}.  
$$

% Step 1.  

Let us set 
$$
\Phi = F_{1}^{(12)}(\Delta\otimes 1)(F_{1})
\left( F_{1}^{(23)}(1\otimes \Delta)(F_{1})\right)^{-1} , 
$$
we have clearly $\Phi \in A^{<0} \otimes A \otimes A^{\geq 0}$. 
Since we also have 
$$
\Phi = \left( (\bar\Delta\otimes 1)(F_{2}) F_{2}^{(12)} \right)^{-1}
(1\otimes \bar\Delta)(F_{2}) F_{2}^{(23)} ,
$$
we also see that $\Phi\in A^{\geq 0} \otimes A \otimes A^{<0}$. 

% Step 2. 

Therefore 
$\Phi = 1 \otimes a \otimes 1$, for a certain $a\in A$. On the 
other hand, as $\Phi$ is obtained by twisting a quasi-Hopf structure, it 
should satisfy the compatibility condition (see \cite{D-qH})
$$
(\Delta_{1} \otimes id \otimes id)(\Phi)
(id \otimes id \otimes \Delta_{1})(\Phi)
=
(\Phi\otimes 1)(id \otimes \Delta_{1}\otimes id)(\Phi) 
(1\otimes \Phi) ,
$$
where $\Delta_{1} = \Ad(F_{1}) \circ \Delta$.  
This implies that 
$$
1\otimes a \otimes a \otimes 1 = (1\otimes a\otimes 1 \otimes 1)
(1\otimes \Delta_{1}(a) \otimes 1) (1\otimes 1 \otimes a \otimes 1), 
$$
and so $\Delta_{1}(a)=1$; applying the counit to one of the factors of
the tensor product where this equality takes place, we obtain $a=1$. 
\hfill \qed  \medskip 

\begin{remark}
An other way to show that $\Phi$ is scalar is the following. 
We can use the third expression of $\Phi$ in \cite{Enr-Rub} Prop. 7.4 to 
show that $\Phi$ belongs to $A \otimes A^{\geq 0}\otimes A$. 
By writing a similar expression for $\Phi$, we get that 
$\Phi\in A \otimes A^{<0} \otimes A$. Together with the fact that $\Phi$
belongs to $1\otimes A \otimes 1$, this shows that $\Phi$ is scalar. 
\end{remark}

\begin{remark} First order computations lead us to believe that $F_{1}$,
resp. $F_{2}$ can be expressed polynomially in terms of the 
$\res_{0}(e^{< 0}(z) \otimes f^{\ge 0}(z))^{n} dz$, resp. of the
$\res_{0}(e^{\ge 0}(z) $ %%%%%%%%%%% LIAISON
$\otimes f^{< 0}(z))^{n} dz$. 
Product formulas for
$F_{1,2}$ can be found in \cite{KhT}. 
\end{remark}

\section{Yangian $RLL$ relations} \label{sect:RLL}

Its follows from Thm. \ref{cocycle} that we can twist the Hopf algebra 
structure $(A,\Delta)$ by $F_{1}$, and get another Hopf algebra structure.
The twisted coproduct is $\Delta_{1} = \Ad(F_{1}) \circ\Delta$. 

Let 
\begin{equation} \label{R:rat}
\cR = q^{D \otimes K} q^{{1\over 2}\sum_{i\ge 0} h_{i} \otimes h_{-i-1}}
q^{\sum_{i\in\ZZ}e_{i}\otimes f_{-i-1}}; 
\end{equation}
this is the universal $R$-matrix for $(A,\Delta)$ (see \cite{Enr-Rub}). 
The universal $R$-matrix for the twisted Hopf algebra $(A, \Delta_{1})$ 
is then $\cR_{1} = F_{1}^{(21)}\cR F_{1}^{-1}$. We then
have the Yang-Baxter equation 
\begin{equation} \label{YB}
\cR_{1}^{(12)}\cR_{1}^{(13)}\cR_{1}^{(23)}
=
\cR_{1}^{(23)}\cR_{1}^{(13)}\cR_{1}^{(12)}. 
\end{equation} 

Recall now the formulas for $2$-dimensional representations of $A$ (see
e.g. \cite{CP}). Let $\zeta$ be a formal variable, $k_{\zeta}$ the field
of formal Laurent power series $\CC((\zeta))$, $\pa_{\zeta}$ the
derivation of $k_{\zeta}$ defined as $d/d\zeta$, and
$k_{\zeta}[\pa_{\zeta}]$ the associated ring of differential operators. 
\begin{lemma}
There is a morphism of algebras $\pi_{\zeta}$ from $A$ to $\End(\CC^{2})
\otimes k_{\zeta}[\pa_{\zeta}][[\hbar]]$, defined by 
$$
\pi_{\zeta}(K) = 0, \quad \pi_{\zeta}(D) = \Id_{\CC^{2}} \otimes
\pa_{\zeta}, 
$$
$$
\pi_{\zeta}(h_{n}) = \pmatrix \left( {2\over{1+q^{\pa_{z}}}}z^{n}\right) 
(\zeta) & 0 \\ 0 & 
- \left( {2\over{1+q^{-\pa_{z}}}}z^{n}\right) (\zeta) \endpmatrix, \quad
n\ge 0,
$$
$$
\pi_{\zeta}(h_{n}) = \pmatrix \left(
{{1-q^{-\pa_{z}}}\over{\hbar\pa_{z}}}z^{n}
\right) 
(\zeta) & 0 \\ 0 & 
- \left( {{q^{\pa_{z}}-1}\over{\hbar\pa_{z}}} z^{n}\right) (\zeta)
\endpmatrix, \quad n <0,
$$
$$
\pi_{\zeta}(e_{n}) = \pmatrix 0 & \zeta^{n} \\ 0 & 0 \endpmatrix, 
\quad
\pi_{\zeta}(f_{n}) = \pmatrix 0 & 0 \\ \zeta^{n} & 0 \endpmatrix, \quad
n\in\ZZ.   
$$
\end{lemma}

\begin{lemma}
We have 
$$
(1\otimes \pi_{\zeta})(\cR_{1}) = L^{\ge 0}(\zeta), 
\quad
(1\otimes \pi_{\zeta})(\cR_{1}^{(21)}) = q^{K\pa_{\zeta}} L^{< 0}(\zeta).  
$$
\end{lemma}

{\em Proof. \/} Let us denote by $U_{\hbar}\N_{\pm}^{\ge i}$ the linear
spans in $U_{\hbar}\N_{\pm}$ of products of more than $i$ factors
$e_{k}$, resp. $f_{k}$. Then the various $\Pi^{\pm}_{*,*}$ preserve the 
$U_{\hbar}\N_{\pm}^{\ge i}$. 
The formulas (\ref{voron}) for $F_{1}$ then imply  
that $F_{1}$ belongs to 
$1+\hbar \sum_{i\ge 0}e_{-i-1}\otimes f_{i} + U_{\hbar}\N_{+}^{\ge 2}
\otimes U_{\hbar}\N_{-}^{\ge 2}$.  
The lemma now follows from the decomposition (\ref{R:rat}), and from the
fact that the $U_{\hbar}\N_{\pm}^{\ge 2}$
are contained in the kernel of $\pi_{\zeta}$. 
\hfill \qed \medskip

\begin{lemma} The image of $\cR_{1}$ by $\pi_{\zeta} \otimes
\pi_{\zeta'}$ is 
$$
(\pi_{\zeta} \otimes \pi_{\zeta'})(\cR_{1}) = A(\zeta,\zeta') 
R^{<0}(\zeta-\zeta'), 
$$
where 
$$
R^{<0}(z) = {1\over{z-\hbar}}\left( z \Id_{\CC^{2}\otimes \CC^{2}} - \hbar
P\right), 
$$
wihere $P$ is the permutation operator of the two factors of
$(\CC^{2})^{\otimes 2}$, and $A(\zeta,\zeta')$ is the formal series 
$\exp\left(\sum_{i\ge 0} ({1\over \pa}{{q^{\pa} -
1}\over{q^{\pa}+1}} \zeta^{i}) \zeta^{\prime -i-1}\right)$. 
\end{lemma}

{\em Proof.} Since the image by $\pi_{\zeta}$ and $\pi_{\zeta'}$
of $U_{\hbar}\N_{\pm}^{\ge 2}$ is equal to zero, and using again the
fact that $F_{1}$ belongs to 
$1+\hbar \sum_{i\ge 0}e_{-i-1}\otimes f_{i} + U_{\hbar}\N_{+}^{\ge 2}
\otimes U_{\hbar}\N_{-}^{\ge 2}$, we find that this image is the same as
that of 
$$
\left( 1+\hbar \sum_{i\ge 0}f_{i} \otimes e_{-i-1}\right) 
q^{{1\over 2} \sum_{i\ge 0} h_{i} \otimes h_{-i-1}}
\left(1-\hbar \sum_{i\ge 0}e_{-i-1}\otimes f_{i} \right). 
$$
Let us denote by $E_{ij}$ the endomorphism of $\CC^{2}$ such that
$E_{ij}v_{\al} = \delta_{\al j}v_{i}$, where $(v_{1},v_{-1})$ is the
standdard basis of $\CC^{2}$. 
We find that  
\begin{align*}
& (\pi_{\zeta} \otimes \pi_{\zeta'})(\cR_{1}) = 
A(\zeta,\zeta') \left( 1 + {\hbar\over{\zeta'-\zeta}} E_{-1,1} \otimes
E_{1,-1} \right) 
\left( E_{1,1} \otimes E_{1,1} + E_{-1,-1} \otimes E_{-1,-1} 
\right. \\ & \left. +
{{\zeta'-\zeta}\over{\zeta'-\zeta+\hbar}} E_{1,1} \otimes
E_{-1,-1}
+ {{\zeta'-\zeta-\hbar}\over{\zeta'-\zeta}} E_{-1,-1} \otimes E_{1,1}
\right)
\left( 1 - {\hbar\over {\zeta-\zeta'}} E_{1,-1} \otimes E_{-1,1}\right); 
\end{align*}
the lemma follows. 
\hfill \qed \medskip 

Define $R^{\ge 0}(z)$ as the inverse of $R^{<0}(z)$. We have 
$$
R^{\ge 0}(z) = {1\over{z+\hbar}}\left( z \Id_{\CC^{2}\otimes \CC^{2}} 
+ \hbar P\right). 
$$

Let us now apply to (\ref{YB}) $1 \otimes \pi_{\zeta} \otimes 
\pi_{\zeta'}$, 
$\pi_{\zeta} \otimes \pi_{\zeta'} \otimes 1$ and $\pi_{\zeta}
\otimes \pi_{\zeta'} \otimes 1$. We find the following relations between
matrices $L^{\pm}(\zeta)$: 

\begin{prop} We have 
\begin{equation} \label{RLL1}
R^{\eta}(\zeta-\zeta')L^{\eta(1)}(\zeta)L^{\eta(2)}(\zeta')
= L^{\eta(2)}(\zeta')L^{\eta(1)}(\zeta) R^{\eta}(\zeta-\zeta') 
\end{equation}
\begin{align} \label{RLL2}
L^{<0(1)}(\zeta) 
R^{<0}(\zeta-\zeta') 
L^{\ge 0(2)}(\zeta')= 
L^{\ge 0(2)}(\zeta') R^{<0}(\zeta-\zeta'-\hbar K) 
L^{<0(1)}(\zeta) {{A(\zeta,\zeta'-K\gamma)}\over{A(\zeta,\zeta')}}, 
\end{align}
$\eta\in\{\ge 0, <0 \}$. 
\end{prop}

\begin{remark}
After analytic continuation in the variables $\zeta,\zeta'$, we see
that $A(\zeta,\zeta')$ only depends on $\zeta - \zeta'$. If we set
$A(\zeta,\zeta') = A(\zeta-\zeta')$, we then have 
$$
A(z) A(z+\hbar) = {z\over{z+\hbar}} , 
$$
so that $A$ is equal to  
$$ 
A(z)= { {\Gamma\left({z\over2\hbar}+{1\over2}\right)^{2}} 
\over 
{ \Gamma\left({z\over{2\hbar}}+1\right)
\Gamma\left({z\over{2\hbar}}\right)} } .
$$ 
\end{remark}

\section{$F_{1}$ and the Yangian coproduct} \label{conjug}

Since $\Delta(A^{\geq 0}) \subset A\otimes A^{\geq 0}$ and $F_{1}$ and
$F_{1}^{-1}$ belong to $A\otimes
A^{\geq 0}$, $\Delta_{1}(A^{\geq 0})\subset A\otimes A^{\geq 0}$. 
On the other hand,
$\Delta_{1} = \Ad(F_{2}^{-1})\circ \bar\Delta$; since 
$\bar\Delta(A^{\geq 0})
\subset A^{\geq 0}\otimes A$, and $F_{2}$ and $F_{2}^{-1}$ belong to 
$A^{\geq 0}\otimes A$, $\Delta_{1}(A^{\geq 0})\subset A^{\geq 0}\otimes A$. 
This shows that $A^{\geq 0}$ is a Hopf subalgebra of $(A,\Delta_{1})$. 

We can show in the same way that $A^{<0}$ is a Hopf subalgebra of
$(A,\Delta_{1})$.  

Therefore it is natural to expect that $\Delta_{1}$ coincides with the 
Yangian coproduct $\Delta_{Yg}$. In this section, we show that this is
indeed the case. 

Since $(A,\Delta_{1}, R_{1})$ is a quasi-triangular Hopf algebra, we
have 
\begin{equation} \label{Delta:R}
(\Delta_{1} \otimes 1)(\cR_{1}) = \cR_{1}^{(13)}\cR_{1}^{(23)}, 
\quad
(1\otimes \Delta_{1})(\cR_{1}) = \cR_{1}^{(13)}\cR_{1}^{(12)}.  
\end{equation}

Apply now $id\otimes id \otimes \pi_{\zeta}$ to the first equation of
(\ref{Delta:R}) and $\pi_{\zeta}\otimes id \otimes id$ to the second
one. We find 
$$
(\Delta_{1}\otimes 1)L^{\ge 0}(\zeta)^{(12)} = 
L^{\ge 0}(\zeta)^{(13)}L^{\ge 0}(\zeta)^{(23)} , 
\quad
(\Delta_{1}\otimes 1)\cL^{< 0}(\zeta)^{(12)} = 
\cL^{< 0}(\zeta)^{(13)}\cL^{< 0}(\zeta)^{(23)} , 
$$
where $\cL^{<0}(\zeta) = q^{K\pa_{\zeta}}L^{<0}(\zeta)$; the last
equation implies,
since $\Delta(K) = K \otimes 1 + 1 \otimes K$, that 
$$
(\Delta_{1}\otimes 1)L^{< 0}(\zeta)^{(12)} = 
L^{< 0}(\zeta-\hbar K_{1})^{(13)} L^{< 0}(\zeta)^{(23)} , 
$$
Since we also have $[D \otimes 1 + 1 \otimes D, F_{1}] = 0$ (the
algebras $A_{\pm}$ being $\ad(D)$-invariant), and $[K\otimes 1 + 1
\otimes K,F_{1}]=0$, and comparing the above formulas with
(\ref{Delta:Yg:1}), (\ref{Delta:Yg:2}), 
we conclude: 

\begin{prop}
$\Delta_{1} = \Ad(F_{1}) \circ \Delta$ coincides with $\Delta_{Yg}$. 
\end{prop}


\frenchspacing
\begin{thebibliography}{10}
{\small
\bibitem{CP} V. Chari, A. Pressley, Quantum affine algebras,
Comm. Math. Phys. 142 (1991), 261-83.    

\bibitem{DF} J. Ding, I.B. Frenkel, Isomorphism of two
realizations of quantum affine algebras $U_{q}(\hat{\frak{gl}}_{n})$,
Comm. Math. Phys. 156 (1993), 277-300.

\bibitem{D-new} V.G. Drinfeld, A new realization of Yangians and
quantized affine algebras, Sov. Math. Dokl. 36 (1988).

\bibitem{D-qH} V.G. Drinfeld, Quasi-Hopf algebras, Leningrad
Math. J. 1:6 (1990), 1419-57. 

\bibitem{toap} B. Enriquez, G. Felder, Elliptic quantum groups 
$E_{\tau,\eta}({\frak{sl}}_{2})$ and quasi-Hopf algebras, in
preparation. 

\bibitem{Enr-Rub} B. Enriquez, V.N. Rubtsov, Quasi-Hopf algebras 
associated with $\frak{sl}_{2}$ and complex curves, preprint
Ecole Polytechnique, no. 1145, q-alg/9608005. 

\bibitem{FRT} L. Faddeev, N. Reshetikhin, L. Takhtajan, Quantization of
Lie groups and Lie algebras, Leningrad Math. J. 1 (1989), 178-201.

\bibitem{Kh}
S. Khoroshkin, Central extension of the Yangian double, preprint
q-alg/9602031. 

\bibitem{KhT} S. Khoroshkin, V. Tolstoy, Yangian double,
Lett. Math. Phys., to appear. 

\bibitem{RS} N. Reshetikhin, M. Semenov-Tian-Shansky, Central extensions
of quantum current groups, Lett. Math. Phys. 19 (1990), 133-42.  

\bibitem{Sm}
F.A. Smirnov, Dynamical symmetries of massive integrable models, 
Jour. Mod. Phys. A7, suppl. 1B (1992), 813-38.
}
\end{thebibliography}
\end{document}